\documentclass{amsart}

\usepackage{amsmath,amsfonts,amssymb,amsthm,amscd,latexsym,euscript}

\newcommand{\dgrm}[1]{\ensuremath{\smash{\underset{\widetilde{\hphantom{#1}}}{#1}} \mathstrut}}

\newcommand{\range}[1]{\ensuremath{\mathrm{codom}({#1})}}

\newcommand{\cal}{\ensuremath{\mathcal}}

\newtheorem {theorem1}{Theorem}[section]
\newtheorem {theorem}[theorem1]{Theorem}

\newtheorem {proposition}[theorem1]{Proposition}
\newtheorem {lemma}[theorem1]{Lemma}
\theoremstyle{definition}
\newtheorem {definition}[theorem1]{Definition}
\theoremstyle{remark}

\newtheorem {question}[theorem1]{Question}

\newcommand{\map}{\ensuremath{\textit{map}}}

\DeclareMathOperator{\colim}{\textup{colim}}

\def\mor{\ensuremath{\textsl{mor}}}
\def\id{\ensuremath{\textsl{id}}}
\def\Lan{\ensuremath{\textup{Lan}}}

\newcommand{\rarrow}{\rightarrow}

\begin{document}
\title[A non-cofibrantly generated model category]
        {An example of a non-cofibrantly generated model category}
\author{Boris Chorny}
\address{Centre de Recerca Matem\`atica, Apartat 50, E-08193 Bellaterra, Spain}
\email{cboris@crm.es}
\thanks{The author is a fellow of the Marie Curie Training Site hosted by the Centre de Recerca
         Matem\`atica (Barcelona), grant nr. HPMT-CT-2000-00075 of the European Commission.}
\date{\today}
\dedicatory{}
\commby{}
%%% ----------------------------------------------------------------------
\begin{abstract}
We show that the model category of diagrams of spaces generated by a
proper class of orbits is not cofibrantly generated. In particular the
category of maps between spaces may be given a non-cofibrantly generated
model structure.
\end{abstract}
%%% ----------------------------------------------------------------------
\maketitle
%%% ----------------------------------------------------------------------
\section{Introduction and formulation of results}
Several examples of non-cofibrantly generated model categories have
appeared recently (see \cite{AHRT}, \cite{ChristHov}, \cite{Isaksen}) in
response to a question stated by Mark Hovey on his home page. In this note
we introduce another family of such examples.
\par
By the \emph{category of spaces}, denoted by ${\cal S}$, we mean the
category of simplicial sets (or compactly generated topological spaces).
There are plenty of model structures on categories of diagrams of spaces,
with different notions of weak equivalences. Some of them are cofibrantly
generated, e.g. for weak equivalences and fibrations being objectwise and
cofibrations obtained by the left lifting property with respect to
trivial fibrations, the corresponding model category is cofibrantly
generated.
\par
Let us remind (from \cite{DF}, \cite{Farjoun}, \cite{DZ}) that a diagram
$O$ of spaces is called an \emph{orbit} if $\colim O = \ast$. The weak
equivalences which we would like to consider arise naturally from the
relation of equivariant homotopy. By the generalized Bredon theorem \cite
{DZ} a map $f:\dgrm X \rarrow \dgrm Y$ is an equivariant homotopy
equivalence between diagrams which are both cofibrant and fibrant iff
$\map(O,f)\!: \map(O,\dgrm X) \rarrow \map(O,\dgrm Y)$ is a weak
equivalence of spaces for any orbit $O$. A model category,
\emph{generated by the collection of orbits}, on diagrams of spaces was
constructed in \cite{Farjoun} with a map $f$ being a weak equivalence
(reps. fibration) iff $\map(O,f)$ is a weak equivalence (resp. fibration)
for any orbit $O$. In the sequel we consider only this model category on
diagrams of spaces. The simplest example of a non-cofibrantly generated
model category is given by the following
\begin{theorem}
\label{main} If $J=(\bullet\rarrow\bullet)$ is the category with two
objects and only one non-identity morphism, then the functor category
${\mathcal M} = {\mathcal S}^J$ of maps of spaces with the model
structure as above is not cofibrantly generated.
\end{theorem}
However, not every small category gives rise to a non-cofibrantly
generated model category of diagrams. For example, if we take $G$ to be a
group, then the above model structure on ${\mathcal S}^G$ is cofibrantly
generated. We conclude the paper by using this example to produce many
other examples of the same nature.
\paragraph{\emph{Acknowledgments}.}
I would like to thank C.~Casacuberta and E.~Dror Farjoun for helpful
conversations about the subject matter of this paper.

\section{Preliminaries}
By an \emph{orbit over a point} in the colimit of a diagram \dgrm X we
mean the pull back of the canonical map $f\!:\dgrm X\rightarrow
\colim\dgrm X$ over $g\!:\ast\rightarrow \colim\dgrm X$. Let $D$  be any
small category enriched over ${\cal S}$. We denote by ${\cal O}$ the
collection of all orbits of $D$. By collection we mean a set or a proper
class with respect to some fixed universe $\mathfrak U$. The operator
$\range{\cdot}$ applied to a collection of maps returns the collection of
ranges. Given a set $I$ of maps in ${\cal M} = {\cal S}^D$, we denote by
$I$-cell the collection of relative $I$-cellular complexes and by
abs-$I$-cell the collection of (absolute) $I$-cellular complexes. See
\cite[2.1.9]{Hovey} for precise definitions.
\begin{definition}
Let ${\cal X} = \{\dgrm X_\alpha\}_{\alpha \in A}$ be a collection of
$D$-shaped diagrams of spaces. The \emph{collection of orbits} of ${\cal
X}$, denoted by $\Omega({\cal X})\subset {\cal O}$, consists of all orbits
$O\in {\cal O}$ such that there exists $\alpha \in A$ and a point $x\in
\colim \dgrm X_\alpha$ with $O$ being the orbit over $x$.
\end{definition}
\begin{lemma} \label{lemma}
Let $I$ be a set of cofibrations in the model category ${\cal M}$ of
$D$-shaped diagrams of spaces. Then $\Omega($abs-$I$-cell$)\subset \Omega
(\range{I})$.
\end{lemma}
\proof Let $\dgrm X \in {\cal M}$ be any $I$-cellular complex. We proceed
by transfinite induction on the $I$-cellular filtration of \dgrm X.
$\dgrm X_{-1} = \emptyset$, hence $\dgrm X_0 \in \range{I}$ and in
particular $\Omega(\dgrm X_0) \subset \Omega(\range{I})$.
\par
Suppose $\dgrm X_\beta$ satisfies $\Omega(\dgrm X_\beta) \subset
\Omega(\range{I})$. We need to show that $\dgrm X_{\beta+1}$, which is
obtained from $\dgrm X_\beta$ by attaching a map $I \ni f:\dgrm A
\hookrightarrow \dgrm B$, satisfies $\Omega(\dgrm X_{\beta+1}) \subset
\Omega(\range{I})$.
\[
\begin{CD}
    {\dgrm A}    @>{\varphi}>>      {\dgrm X_{\beta}}\\
        @VfV\textup{\quad\tiny push-out}V                             @VV{f'}V\\
    {\dgrm B}    @>>>             {\dgrm X_{\beta+1}}
\end{CD}
\]
Let $O_s$ be an orbit over a point $s \in \colim \dgrm X_{\beta+1} =
\colim\dgrm X_\beta \amalg_{\colim\dgrm A}\colim\dgrm B$. Considering two
cases, $s\in \colim\dgrm X_\beta\subset\colim\dgrm X_{\beta+1}$ and
$s\notin \colim\dgrm X_\beta$, we find out that in the first case $O_s$
equals the corresponding orbit of $\dgrm X_{\beta}$ and in the second
case $O_s$ is some orbit of \dgrm B. This follows immediately from the
fact that the diagrams
\[
\begin{CD}
    {\dgrm X_\beta}         @>{f'}>>    {\dgrm X_{\beta+1}}\\
        @VVV                            @VVV\\
    {\colim \dgrm X_\beta}  @>>>    {\colim \dgrm X_{\beta+1}}
\end{CD}
\qquad
\begin{CD}
    {\dgrm B/\dgrm A}             @>\cong>>    {\dgrm X_{\beta+1}/\dgrm X_\beta}\\
        @VVV                                               @VVV\\
    {\colim(\dgrm B/\dgrm A)}     @>\cong>>    {\colim(\dgrm X_{\beta+1}/\dgrm X_\beta)}
\end{CD}
\]
are pull-backs. The first square is a pull-back by \cite[2.1]{Farjoun} and
the second by the observation that horizontal maps are isomorphisms.
Hence $\Omega(\dgrm X_{\beta+1}) \subset \Omega (\range{I})$.
\par
Obviously, if $\beta$ is a limit ordinal, then
\[
\Omega(\dgrm X_\beta) = \bigcup_{\lambda < \beta} \Omega(\dgrm X_\lambda)
\subset \Omega (\range {I}).
\]
Hence $\Omega($abs-$I$-cell$) \subset \Omega (\range{I})$. \qed
\section{Proof of Theorem \ref{main}}
Let us prove first a slightly more general result.
\begin{proposition} \label{general}
Let $D$ be a small category enriched over $\cal S$ which admits a proper
class of orbits ${\cal O_D}$. Then the model category ${\cal M}$ on the
$D$-shaped diagrams of spaces generated by the orbits is not cofibrantly
generated.
\end{proposition}
\begin{proof}
We argue by contradiction. Suppose the model category $\cal M$ generated
by the proper class of orbits ${\cal O_D}$ is cofibrantly generated. Let
$I$ be the set of generating cofibrations, then any cofibration is a
retract of an $I$-cellular map. (This follows from Quillen's small object
argument; see \cite[2.1.15]{Hovey}.) In particular, any orbit of $\cal
O_D$ is a retract of an $I$-cellular complex. Hence any orbit is a
retract of some orbit of an $I$-cellular space. But by \ref{lemma} the
whole collection of orbits of $I$-cellular complexes form a set, hence
the contradiction.
\end{proof}
\par
In particular, the model category ${\cal S}^J$ is generated by the proper
class of orbits ${\cal O}_J =\{X\rarrow \ast\}$, where $X$ runs through
all the objects of $\cal S$. Therefore, by Proposition \ref{general},
${\cal S}^J$ is not cofibrantly generated, hence the main result
\ref{main}.
\section{More examples}
Let us conclude by giving more examples of non-cofibrantly generated
model categories. Proposition \ref{general} implies that $\cal M = {\cal S
}^D$ is not cofibrantly generated iff ${\cal O}_D$ is a proper class. We
have already indicated in the introduction that if we take $D=G$ to be a
group, then ${\cal O}_G = \{G/H | H<G\}$ is a set, hence ${\cal S}^G$ is
cofibrantly generated. The same holds for groupoids. However, the
following proposition provides us with a large family of examples.
\begin{proposition}
Let $D$ be a small category which admits a fully faithful functor $i:K
\rarrow D$, where $K$ is a category with two objects $k_1, k_2$, at least
one arrow $f:k_1\rarrow k_2$ and no arrows in the opposite direction. Then
${\cal O}_D$ consists of a proper class of orbits.
\end{proposition}
\begin{proof}
First define for each space $X\in {\cal S}$ an orbit over $K$, i.e. a
functor $T_X\!:K\rarrow {\cal S}$, by $T_X(k_1) = X,\; T_X(k_2) = \ast,\;
T_X(g) = \id_X$ for any $g\in \mor(k_1, k_1)$ and $T_X$ on the elements
of $\mor(k_1, k_2),\; \mor(k_2, k_2)$ has a unique definition, since
$\ast$ is the final object of ${\cal S}$. Obviously $T_X$ is an orbit
over $K$.
\par
Next we define for each $T_X$ a $D$-orbit $O_X$ by extending the
definition of $T_X$ to the whole $D$. More precisely, $O_X = \Lan_i\,
T_X$. We need to check that $O_X$ is an orbit. It follows from the fact
that colimit is itself a left Kan extension along a functor into the
trivial category. But any two left Kan extensions commute since they may
be represented as coends, and for the coends there is a ``Fubini''
theorem. See \cite[X]{MacLane} for the details.
\par
The functor $i$ is taken to be fully faithful, hence $O_X(i(k_1)) = X$;
therefore we obtain a proper class of $D$-orbits of the form $O_X$.
\end{proof}
\begin{question}
Let $D$ be a monoid which is not a group. Is ${\cal S}^D$ cofibrantly
generated?
\end{question}

\end{document}